\documentclass[12pt,twoside]{article}
\usepackage{amssymb}
\usepackage{amsmath}

\date{}

\begin{document}
\centerline {\Large{  A Characterization of Inner Product Spaces }
 Related to}

\centerline{}

\centerline {\Large{ the Skew-Angular Distance}}

\centerline{}
\centerline{ Hossein Dehghan
}
\centerline{}
\centerline{Department of Mathematics, Institute for
Advanced Studies in Basic
 Sciences
 }
\centerline{ (IASBS), Gava Zang,  Zanjan 45137-66731, Iran}
\centerline{ Email: h$_{-}$dehghan@iasbs.ac.ir, hossein.dehgan@gmail.com}
\centerline{}

\newtheorem{theorem}{Theorem}[section]
\newtheorem{lemma}[theorem]{Lemma}
\newtheorem{proposition}[theorem]{Proposition}
\newtheorem{corollary}[theorem]{Corollary}
\newtheorem{definition}[theorem]{Definition}
\newtheorem{example}[theorem]{Example}
\newtheorem{xca}[theorem]{Exercise}
\newtheorem{remark}[theorem]{Remark}
\newtheorem{algorithm}{Algorithm}
\numberwithin{equation}{section}

\begin{abstract}
 \par
  A new refinement of the triangle inequality is presented in normed linear spaces. Moreover, a simple characterization of inner product spaces is obtained by using the skew-angular distance.

\end{abstract}

{\bf Keywords:} Triangle inequality, characterization of inner product spaces, angular distance

\section{ Introduction}
\par
In 2006, Maligranda \cite[Theorem 1]{Mal1} (see also \cite{Mal2}) introduced the following strengthening of the triangle inequality and its reverse: For any nonzero vectors $x$ and $y$ in a real normed linear space $X = (X, \|.\|)$ it is true that
\begin{eqnarray}\label{Ma re tri up}
\|x+y\|\leq\|x\|+\|y\|-\left(2-\left\|\frac{x}{\|x\|}+\frac{y}{\|y\|}\right\|\right)
\min\{\|x\|,\|y\|\}
\end{eqnarray}
and
\begin{eqnarray}\label{Ma re tri low}
\|x+y\|\geq\|x\|+\|y\|-\left(2-\left\|\frac{x}{\|x\|}+\frac{y}{\|y\|}\right\|\right)
\max\{\|x\|,\|y\|\}.
\end{eqnarray}
Also, the author used  (\ref{Ma re tri up}) and (\ref{Ma re tri low}) for the following estimation of the \emph{angular
distance}  $\alpha[x,y]=\|\frac{x}{\|x\|}-\frac{y}{\|y\|}\|$ between two nonzero elements $x$ and $y$ in $X$ which was defined by Clarkson  in \cite{Clarc}.
\begin{eqnarray}\label{angu}
\frac{\|x-y\|-|\ \|x\|-\|y\|\ |}{\min\{\|x\|,\|y\|\}}\leq\alpha[x,y]\leq\frac{\|x-y\|+|\ \|x\|-\|y\|\ |}{\max\{\|x\|,\|y\|\}}.
\end{eqnarray}
The right hand of estimate (\ref{angu}) is a refinement of the Massera-Schaffer inequality proved in 1958
(see \cite[ Lemma 5.1]{MS}): for nonzero
vectors $x$ and $y$ in $X$ we have that $\alpha[x,y]\leq\frac{2\|x-y\|}{\max\{\|x\|,\|y\|\}}$ , which is stronger than
the Dunkl-Williams inequality $\alpha[x,y]\leq\frac{4\|x-y\|}{\|x\|+\|y\|}$ proved in \cite{DW}. In the same paper, Dunkl and Williams proved that the constant 4 can be replaced by 2 if and only if $X$ is an inner product space.
\par
The main aim of this paper is to obtain a new and simple characterization of inner product spaces. To proceed in this direction we first present a refinement of the triangle inequality in normed linear spaces and introduce the notion of skew-angular distance. Next, we compare the angular distance and skew-angular distance with each other.
\section{ A refinement of triangle inequality}
\par
 we start with the following strengthening of the triangle inequality.
\begin{theorem}
For nonzero vectors $x$ and $y$ in a real normed linear space $X = (X, \|.\|)$ we have
\begin{eqnarray}\label{my re tri up}
\|x+y\|\leq\|x\|+\|y\|-\left(\frac{\|x\|}{\|y\|}+\frac{\|y\|}{\|x\|}-\left\|\frac{x}{\|y\|}
+\frac{y}{\|x\|}\right\|\right)
\min\{\|x\|,\|y\|\}
\end{eqnarray}
and
\begin{eqnarray}\label{my re tri low}
\|x+y\|\geq\|x\|+\|y\|-\left(\frac{\|x\|}{\|y\|}+\frac{\|y\|}{\|x\|}-\left\|\frac{x}{\|y\|}
+\frac{y}{\|x\|}\right\|\right)
\max\{\|x\|,\|y\|\}.
\end{eqnarray}
\end{theorem}
\textbf{ Proof.} Without loss of generality we may assume that $\|x\|\leq\|y\|$. Then, by the triangle inequality,
\begin{align}\label{}
\nonumber\|x+y\|&=\left\|\frac{\|x\|}{\|y\|}x+\frac{\|x\|}{\|x\|}y+\left(1-\frac{\|x\|}{\|y\|}
\right)x\right\|\\
\nonumber&\leq\|x\|\left\|\frac{x}{\|y\|}+\frac{y}{\|x\|}\right\|+\|x\|-\frac{\|x\|^2}{\|y\|}
\\
\nonumber&=\|x\|+\|x\|\left(\left\|\frac{x}{\|y\|}+\frac{y}{\|x\|}\right\|
-\frac{\|x\|}{\|y\|}\right)
\\
\nonumber&=\|x\|+\|y\|+\|x\|\left(\left\|\frac{x}{\|y\|}+\frac{y}{\|x\|}\right\|
-\frac{\|x\|}{\|y\|}-\frac{\|y\|}{\|x\|}\right)
\end{align}
which establishes estimate (\ref{my re tri up}). Similarly, the computation
\begin{align}\label{}
\nonumber\|x+y\|&=\left\|\frac{\|y\|}{\|y\|}x+\frac{\|y\|}{\|x\|}y+\left(1-\frac{\|y\|}{\|x\|}
\right)y\right\|\\
\nonumber&\geq\|y\|\left\|\frac{x}{\|y\|}+\frac{y}{\|x\|}\right\|-\left(
\frac{\|y\|^2}{\|x\|}-\|y\|\right)\\
\nonumber&=\|y\|+\|y\|\left(\left\|\frac{x}{\|y\|}+\frac{y}{\|x\|}\right\|-\frac{\|y\|}{\|x\|}\right)
\\
\nonumber&=\|x\|+\|y\|+\|y\|\left(\left\|\frac{x}{\|y\|}+\frac{y}{\|x\|}\right\|
-\frac{\|x\|}{\|y\|}-\frac{\|y\|}{\|x\|}\right)
\end{align}
gives inequality (\ref{my re tri low}). $\Box$
\par
The following examples show that neither our refinement nor Maligranda's refinement of the triangle inequality is always better.
\begin{example}  Let $X$  be the normed space $\mathbb{R}$ with the norm  $\|x\|=|x|$. Then, for $x=1$ and $y=-2$ we have
\begin{eqnarray*}
2-\left\|\frac{x}{\|x\|}+\frac{y}{\|y\|}\right\|=2>1
=
\frac{\|x\|}{\|y\|}+\frac{\|y\|}{\|x\|}-\left\|\frac{x}{\|y\|}+\frac{y}{\|x\|}\right\|.
\end{eqnarray*}
\end{example}
\begin{example}  Let $X=\mathbb{R}^2$ with the norm of $x=(a,b)$ given by $\|x\|=|a|+|b|$. Take $x=(3/4,3/4)$ and $y=(-1,0)$, then
 $\|x\|=3/2$, $\|y\|=1$. Therefore,
$$\frac{x}{\|x\|}=(\frac{1}{2},\frac{1}{2}),\ \ \ \ \frac{y}{\|y\|}=(-1,0),\ \ \ \frac{x}{\|y\|}=(\frac{3}{4},\frac{3}{4}),\ \ \ \ \frac{y}{\|x\|}=(-\frac{2}{3},0)$$
and
\begin{eqnarray*}
2-\left\|\frac{x}{\|x\|}+\frac{y}{\|y\|}\right\|=1<\frac{4}{3}
=
\frac{\|x\|}{\|y\|}+\frac{\|y\|}{\|x\|}-\left\|\frac{x}{\|y\|}+\frac{y}{\|x\|}\right\|.
\end{eqnarray*}
\end{example}
\par
We can gather estimates (\ref{my re tri up}) and (\ref{my re tri low}) together as
\begin{align}
\nonumber\|x+y\|&+\left(\frac{\|x\|}{\|y\|}+\frac{\|y\|}{\|x\|}-\left\|\frac{x}{\|y\|}
+\frac{y}{\|x\|}\right\|\right)\min\{\|x\|,\|y\|\}\leq\|x\|+\|y\|\\
\nonumber&\leq\|x+y\|+\left(\frac{\|x\|}{\|y\|}+\frac{\|y\|}{\|x\|}-\left\|\frac{x}{\|y\|}
+\frac{y}{\|x\|}\right\|\right)\max\{\|x\|,\|y\|\}.
\end{align}
Also, we use them as the estimates for a distance in normed linear spaces which we call \emph{skew-angular distance}.
\begin{definition}
 For two nonzero elements $x$ and $y$ in a real normed linear space $X=(X,\|.\|)$  the distance
\begin{eqnarray}
\beta[x,y]=\left\|\frac{x}{\|y\|}-\frac{y}{\|x\|}\right\|
\end{eqnarray}
is called skew-angular distance between  $x$ and $y$.
\end{definition}
\begin{corollary}
For any nonzero elements $x$ and $y$ in a real normed linear space $X=(X,\|.\|)$ we have
\begin{eqnarray}\label{sk up}
\beta[x,y]\leq \frac{\|x-y\|}{\max\{\|x\|,\|y\|\}}+\frac{|\ \|x\|-\|y\|\ |}{\min\{\|x\|,\|y\|\}}
\end{eqnarray}
and
\begin{eqnarray}\label{sk up}
\beta[x,y]\geq \frac{\|x-y\|}{\min\{\|x\|,\|y\|\}}-\frac{|\ \|x\|-\|y\|\ |}{\max\{\|x\|,\|y\|\}}
\end{eqnarray}
\end{corollary}
\textbf{ Proof.}
  Estimate (\ref{my re tri low}) implies that
\begin{align}
\nonumber\left\|\frac{x}{\|y\|}
-\frac{y}{\|x\|}\right\|\max\{\|x\|,\|y\|\}\leq\|x-y\|-\|x\|-\|y\|+\left(\frac{\|x\|}{\|y\|}
+\frac{\|y\|}{\|x\|}\right)\max\{\|x\|,\|y\|\}.
\end{align}
Without loss of generality we may assume that $\|x\|\leq\|y\|$. Then
\begin{align}
\nonumber\left\|\frac{x}{\|y\|}
-\frac{y}{\|x\|}\right\|\|y\|\leq\|x-y\|+\frac{\|y\|}{\|x\|}(\|y\|-\|x\|)
\end{align}
and so
\begin{align}
\nonumber\left\|\frac{x}{\|y\|}
-\frac{y}{\|x\|}\right\|\leq\frac{\|x-y\|}{\|y\|}+\frac{|\ \|y\|-\|x\|\ |}{\|x\|}.
\end{align}
Similarly, inequality (\ref{my re tri up}) gives that
\begin{align}
\nonumber\left\|\frac{x}{\|y\|}
-\frac{y}{\|x\|}\right\|\|x\|\geq\|x-y\|-\frac{\|x\|}{\|y\|}(\|y\|-\|x\|)
\end{align}
and so
\begin{align}
\nonumber\left\|\frac{x}{\|y\|}
-\frac{y}{\|x\|}\right\|\geq\frac{\|x-y\|}{\|x\|}-\frac{|\ \|y\|-\|x\|\ |}{\|y\|}.
\end{align}
This completes the proof. $\Box$\\
\par
 Estimates (\ref{my re tri up}) and (\ref{my re tri low}) mean for the skew-angular distance that
\begin{eqnarray*}\label{}
\frac{\|x-y\|}{\min\{\|x\|,\|y\|\}}-\frac{|\ \|x\|-\|y\|\ |}{\max\{\|x\|,\|y\|\}}\leq \beta[x,y]\leq \frac{\|x-y\|}{\max\{\|x\|,\|y\|\}}+\frac{|\ \|x\|-\|y\|\ |}{\min\{\|x\|,\|y\|\}}.
\end{eqnarray*}
Since $|\ \|x\|-\|y\|\ |\leq \|x-y\|$, we obtain the estimate
\begin{eqnarray}\label{m-type} \beta[x,y]\leq\left( \frac{1}{\max\{\|x\|,\|y\|\}}+\frac{1}{\min\{\|x\|,\|y\|\}}\right)\|x-y\|=\left( \frac{1}{\|x\|}+\frac{1}{\|y\|}\right)\|x-y\|.
\end{eqnarray}
 The constant $1$ in the estimate (\ref{m-type}) is the best possible even for an inner product space. In fact, consider $X=\mathbb{R}$ with the norm of $x$ given by $\|x\|=|x|$. Take $x=-1$ and $y=\epsilon$, where $\epsilon>0$ is small. Then
 \begin{eqnarray*}\label{}
\beta[x,y]=\epsilon+\frac{1}{\epsilon} \ \ \ and \ \ \ \left( \frac{1}{\|x\|}+\frac{1}{\|y\|}\right)\|x-y\|=\left(1+\frac{1}{\epsilon}\right)(1+\epsilon).
\end{eqnarray*}
Hence
 \begin{eqnarray*}\label{}
 \beta[x,y] \frac{\|x\|\|y\|}{(\|x\|+\|y\|)\|x-y\|}
=\frac{1+\epsilon^2}{(1+\epsilon)^2} \rightarrow 1
\end{eqnarray*}
as $\epsilon\rightarrow 0^+$.
\section{ Characterization of inner product spaces}
In this section we compare the norm-angular distance $\alpha[x,y]$ with the skew-angular distance $\beta[x,y]$.
 The next theorem due to Lorch will be useful in the sequel.
\begin{theorem}\label{Lorch th}
(See \cite{Lorch}.) Let $(X,\|.\|)$ be a real normed linear space. Then the following statements are mutually equivalent:
\begin{enumerate}
  \item[(i)] For each $x,y\in X$ if $\|x\|=\|y\|$, then $\|x+y\|\leq \|\gamma x+\gamma^{-1}y\|$ (for all $\gamma\neq 0$).
  \item[(ii)] For each $x,y\in X$ if $\|x+y\|\leq \|\gamma x+\gamma^{-1}y\|$ (for all $\gamma\neq 0$), then $\|x\|=\|y\|$.
  \item[(iii)] $(X,\|.\|)$  is an inner product space.
  \end{enumerate}
\end{theorem}
\begin{theorem}\label{chara inner}
Let $(X,\|.\|)$ be a real normed linear space. Then
$ (X, \|.\|)$ is  an inner product space if and only if for each nonzero elements $x$ and $y$ in $X$,
\begin{eqnarray}\label{alpha beta}
\alpha[x,y]\leq\beta[x,y].
\end{eqnarray}

\end{theorem}\label{comparison p=0}
\textbf{ Proof.} Let $X=(X,\langle.,.\rangle)$ be inner product space, $x,y\in X$ and $x,y\neq0$. We consider that
\begin{align}\label{}
\nonumber\left\|\frac{x}{\|y\|}-\frac{y}{\|x\|}\right\|^2&-\left\|\frac{x}{\|x\|}-\frac{y}{\|y\|}\right\|^2
=\left\langle\frac{x}{\|y\|}-\frac{y}{\|x\|},\frac{x}{\|y\|}-\frac{y}{\|x\|}\right\rangle-\left\langle\frac{x}{\|x\|}
-\frac{y}{\|y\|},\frac{x}{\|x\|}-\frac{y}{\|y\|}\right\rangle\\
\nonumber&=\frac{\|x\|^2}{\|y\|^2}+\frac{\|y\|^2}{\|x\|^2}-\frac{2\langle x,y\rangle}{\|x\|\|y\|}-\left(\frac{\|x\|^2}{\|x\|^2}
+\frac{\|y\|^2}{\|y\|^2}-\frac{2\langle x,y\rangle}{\|x\|\|y\|}\right)\\
\nonumber&=\frac{\|x\|^2}{\|y\|^2}+\frac{\|y\|^2}{\|x\|^2}-2=\left(\frac{\|x\|}{\|y\|}-\frac{\|y\|}{\|x\|}\right)^2\geq0.
\end{align}
This proves the necessity.
\par
 To prove the sufficiency let $x,y\in X$, $\|x\|=\|y\|$ and $\gamma\neq 0$. From Theorem \ref{Lorch th} it is enough to prove that  $\|x+y\|\leq \|\gamma x+\gamma^{-1}y\|$. If $x=0$ or $y=0$, then the proof is clear. Let $x\neq 0$, $y\neq 0$ and $\gamma>0$. Applying inequality (\ref{alpha beta} ) to  $\gamma^{\frac{1}{2}} x$ and $-\gamma^{-\frac{1}{2}}y$ instead of $x$ and $y$, respectively, we obtain
\begin{eqnarray}\label{}\nonumber
\left\|\frac{\gamma^{\frac{1}{2}}x}{\gamma^{\frac{1}{2}}\|x\|}+\frac{\gamma^{-\frac{1}{2}}y}{\gamma^{-\frac{1}{2}}\|y\|}\right\|
\leq\left\|\frac{\gamma^{\frac{1}{2}}x}{\gamma^{-\frac{1}{2}}\|y\|}+\frac{\gamma^{-\frac{1}{2}}y}{\gamma^{\frac{1}{2}}\|x\|}\right\|.
\end{eqnarray}
Thus
\begin{eqnarray}\label{}\nonumber
\left\|\frac{x}{\|x\|}+\frac{y}{\|y\|}\right\|
\leq\left\|\frac{\gamma x}{\|y\|}+\frac{\gamma^{-1}y}{\|x\|}\right\|.
\end{eqnarray}
Since $\|x\|=\|y\|\neq0$, then
$$\|x+y\|\leq \|\gamma x+\gamma^{-1}y\|.$$
Now, let $\gamma$ be negative. Put $\mu=-\gamma>0$. From the positive case we get
$$\|x+y\|\leq  \|\mu x+\mu^{-1}y\|=\|\gamma x+\gamma^{-1}y\|.$$
This completes the proof. $\Box$\\
\textbf{Acknowledgment}\\
The author thanks Prof. J. Rooin for his valuable suggestions which improved the original manuscript.

\end{document}